\documentclass[12pt]{article}
\usepackage{amssymb,amsmath,amsfonts,amsthm,amsxtra}
\usepackage{bm}
\usepackage{cite}
\usepackage{array}
\usepackage{url}
\usepackage{hhline}

\def\suml {\mathop{\sum}   \limits}                             %+
\def\ee{e}                                                      %+
\def\e{\mathrm e}                                               %+
\def\a{\alpha}                                                  %+
\def\th{\theta}                                                 %+
\def\g{\gamma}                                                  %+
\def\d{\delta}                                                  %+
\def\cA{\Check A}                                               %++
\def\om{{\mathrm r}}                                            %++
\def\omb{{\bar\om}}                                             %++
\def\ka{\kappa}                                                 %+
\def\RR{\mathop{\cal R}\nolimits}                               %+
\def\CC{\mathop{\cal C}\nolimits}                               %+
\def\cdc{,\ldots,}                                              %+
\def\1n{1\cdc n}                                                %+
\def\eq#1{\begin{equation}#1\end{equation}}                     %+
\def\eqs*#1{\begin{eqnarray*}#1\end{eqnarray*}}                 %+
\def\eqss#1{\begin{eqnarray}#1\end{eqnarray}}                   %+
\newtheorem{thm}{Theorem}{\bfseries}{\itshape}                  %+
\newtheorem{corol}{Corollary}{\bfseries}{\itshape}              %+
\newtheorem{defin}{Definition}{\bfseries}{\upshape}             %+
\newtheorem{lemma}{Lemma}{\bfseries}{\itshape}                  %+
{\bfseries}{\itshape}              %+
\def\proof{{\noindent\bf Proof. }}                              %+
\def\R{{\mathbb R}}                                             %+
\def\ovec{\overrightarrow}                                      %+
\def\tot{\rightleftarrows}                                      %+
\def\T{{\xz\rm\scriptscriptstyle T}\xz}                         %+
\def\LW{{\rm\scriptscriptstyle LW}}                             %+
\def\W{{\hspace{.05em}\rm\scriptscriptstyle W}}                 %+
\def\deW{d^{\hspace{.05em}\ee\hspace{-.05em}\rm\scriptscriptstyle W}_\a}
\def\LeW{{\hspace{.04em}{\rm\scriptscriptstyle L}\ee\hspace{-.05em}{\rm\scriptscriptstyle W}}}%+
\def\s{{\hspace{.07em}\mathrm s}}                               %+
\def\ws{{\hspace{.05em}\rm ws}}                                 %+
\def\SS{{\cal S}}                                               %+
\def\xy{\hspace{.07em}}                                         %+
\def\xz{\hspace{-.07em}}                                        %+
\def\ms{\mathstrut}                                             %+
\def\mz{{\mathstrut}}                                           %+
\def\pha{\phantom{\frac{1}{1}}}                                 %+
\def\diag{\operatorname{diag}}                                  %+
\def\_#1{{^{}_{#1}}}                     %+
\newtheorem{remark}{Remark}{\bfseries}{\upshape}                %+
\newtheorem{example}{Example}{\bfseries}{\upshape}              %+
\def\j{{\bar\jmath}}                                            %+
\def\i{{\bar\imath}}                                            %+
\def\jj{{\bar\jmath\bar\jmath}}                                 %+
\def\ii{{\bar\imath\hspace{.075ex}\bar\imath}}                  %+
\def\ji{{\bar\jmath\bar\imath}}                                 %+
\def\jjt{{\widetilde{\jmath\jmath\,}}}                          %+
\def\iit{{\widetilde{\imath\imath}}}                            %+
\def\oY{{\mathaccent 23 Y}\xz}                                  %+
\def\tY{\widetilde{Y}\xz}                                       %+
\def\l{\ell}                                                    %+
\def\A{{\cal A}}                                                %+
\def\dr{d^{\hspace*{.04em}\om}}                                 %+
\def\dLW{d^{\hspace{.05em}\LW}\xz}                              %+
\def\Up#1{\vspace{-#1em}}                                       %+
\def\ti{\widetilde}                                             %+
\def\L{{\cal L}}                                                %+
\def\Ker{\operatorname{Ker}}                                    %+
\def\ak{_\mz}                                                   %+
\def\aj{^\mz}                                                   %+
\def\ajj{^{\mz^\mz}}                                            %+
\def\akk{_{\mz_\mz}}                                            %+
\def\hy{\hspace{.04em}}                                         %+
\def\hz{\hspace{-.04em}}                                        %+
\def\uv{{\hy\overline{\hz uv\hz}\hy}}                           %+
\def\vu{{\hy\overline{\hz vu\hz}\hy}}                           %+
\def\u{{\hy\overline{\hz u\hz}\hy}}                             %+
\def\v{{\hy\overline{\hz v\hz}\hy}}                             %+

\author{\Up{.1}Pavel Chebotarev\\
       {\Up{.3}\normalsize Institute of Control Sciences of the Russian Academy of Sciences}\\%\Up{.3}
       {\Up{.3}\normalsize 65 Profsoyuznaya Street, Moscow 117997, Russia}\\%\Up{.3}
       {\small\tt chv@member.ams.org}
       }

\title{\vspace{-2.8em}The Walk Distances in Graphs\\
      {\vspace{-0.4em}\small\emph{To the memory of Gerald Subak-Sharpe$,$ Holocaust survivor$,$ a favorite professor} (1925--2011)}\Up{.5}}
\date{}

\begin{document}
\def\baselinestretch{0.94}
\maketitle

\vspace{-2em}
\begin{abstract}
The walk distances in graphs are defined as the result of appropriate transformations of the $\sum_{k=0}^\infty(tA)^k$ proximity measures, where $A$ is the weighted adjacency matrix of a graph and $t$ is a sufficiently small positive parameter. The walk distances are graph-geodetic; moreover, they converge to the shortest path distance and to the so-called long walk distance as the parameter $t$ approaches its limiting values. We also show that the logarithmic forest distances which are known to generalize the resistance distance and the shortest path distance are a specific subclass of walk distances. On the other hand, the long walk distance is equal to the resistance distance in a %certain
transformed graph.


\end{abstract}

\section{Introduction}

The classical distances for graph vertices are the shortest path distance~\cite{BuckleyHarary90}, the resistance distance \cite{Sharpe67,Sharpe67a,Sharpe67b,GvishianiGurvich87,KleinRandic93}, which is proportional to the commute time distance~\cite{GobelJagers74}, and the square root version of the resistance distance \cite{MooreSubak-Sharpe68,Subak-Sharpe89,Subak-Sharpe90}. Recently, a need for a wider variety of graph distances has been strongly felt (see, e.g.,~\cite{YenSaerensShimbo08,DezaDeza09,LuxburgRadlHein09,Tang10PhD,Estrada11book} among many others).

Recall the well-known fact that the shortest path distance and the resistance distance coincide on each tree. In particular, for every path, the resistance distance between every two adjacent vertices is one, as well as the shortest path distance. However, in some applications two central adjacent vertices in a path may be considered as {\xz}being \xz\emph{closer}{\xz} to each other than two peripheral adjacent vertices are as there are more walks (of length\:3, 5, etc.) connecting two central vertices. Such a ``gravitational'' property holds for the forest distances~\cite{CheSha00}. In some other applications, a terminal vertex in a path can be considered as being closer to its neighbor than two central adjacent vertices are. For example, if someone has a single friend, then this friendship is often stronger than that between persons having more friends. This heuristic is supported by the logarithmic forest distances~\cite{Che08dist}.

In \cite{Che10bottl}, a general framework was proposed for constructing graph-geodetic distances\footnote{In this paper, a \emph{distance\/} is assumed to satisfy the axioms of metric.} (a~distance $d(i,j)$ for graph vertices is graph-geodetic whenever $d(i,j)+d(j,k)=d(i,k)$ if and only if every path connecting $i$ and $k$ visits~$j$). Namely, it has been shown that if a matrix $S=(s_{ij})$ produces a strictly positive transitional measure on a graph $G$ (i.e., $s_{ij}\,s_{\!jk}\le s_{ik}\,s_{\!jj}$ for all vertices $i$, $j$, and $k,$ while $s_{ij}\,s_{\!jk}=s_{ik}\,s_{\!jj}$ if and only if every path from $i$ to $k$ visits~$j$), then the logarithmic transformation $h_{ij}=\ln s_{ij}$ and the inverse covariance mapping $d_{ij}=h_{ii}+h_{jj}-h_{ij}-h_{ji}$ convert $S$ into the matrix of a graph-geodetic distance. In the case of digraphs, five transitional measures were found in \cite{Che10bottl}, namely,  the ``connection reliability'', the ``path accessibility'' with a sufficiently small parameter, the ``walk accessibility'', and two versions of the ``forest accessibility''. The distances produced by the forest accessibility on weighted multigraphs (networks) were studied in~\cite{Che08dist}.

In \cite{CheSha98} we applied the inverse covariance mapping to the matrices of walk weights $\sum_{k=0}^\infty (tA)^k,$ where $A$ is the adjacency matrix of a graph, and showed that this leads to distances whenever the positive parameter $t$ is sufficiently small. However, these distances~are not graph-geodetic and some of their properties are exotic (see Section\,\ref{s_example}).

In the present paper, we study the class of graph-geodetic \emph{walk distances}, which involves the logarithmic transformation. Sections~\ref{s_notat} and \ref{s_prel} contain definitions and preliminaries, in Section~\ref{s_trans} the walk distances are expressed in terms of commute cycles and via block matrix operations. Sections~\ref{s_short} and \ref{s_long} are devoted to two limiting cases of walk distances: the short walk distance coincides with the classical shortest path distance, while the long walk distance is original. In Section~\ref{s_e-modi}, we consider modified walk distances (the ``$\ee$-walk distances'') which generalize the classical \emph{weighted\/} shortest path distance. In Section~\ref{s_forest}, it is shown that adding ``balancing loops'' %enables one to generate the logarithmic forest distances in the framework
converts the logarithmic forest distances into a subclass of walk distances. This implies, in particular, that the resistance distance is also a limiting walk distance, as shown in Section\:\ref{s_walkresi}. In Section\:\ref{s_example}, several graph metrics are compared on a simple example.

\section{Notation}
\label{s_notat}

In the graph definitions we mainly follow~\cite{Harary69}.
Let $G$ be a weighted multigraph (a weighted graph where multiple edges are allowed) with vertex set $V(G)=V,$ $|V|=n>1,$ %$V(G)=\{\1n\},$ $n>1$
and edge set~$E(G)$. Loops are allowed; throughout the paper we assume that $G$ is connected. For brevity, we call $G$ a \emph{graph}.
For ${i,j\in V(G),}$ let $n_{ij}\in\{0,1,\ldots\}$ be the number of edges incident to both $i$ and $j$ in~$G$; for every ${q\in\{\1n_{ij}\}}$, $w_{ij}^q>0$ is the weight of the $q\/$th edge of this type. Let
\eq{
\label{e_aij}%+
a_{ij}=\sum_{q=1}^{n_{ij}}w_{ij}^q
}
(if $n_{ij}=0,$ we set $a_{ij}=0$) and $A=A(G)=(a_{ij})_{n\times n};$ $A$~is the symmetric \emph{weighted adjacency matrix\/} of~$G$.
In this paper, all matrix entries are indexed by the vertices of~$G.$ This remark is essential when submatrices are considered: say, ``the $i$th column'' of a submatrix of~$A$ means ``the column corresponding to the vertex $i$ of~$G$'' rather than just the ``column number~$i$'', which may differ.

By the \emph{weight\/} of a graph $G$, $w(G)$, we mean the product %\footnote{Note that in some papers, the weight of a graph is defined as the \emph{sum\/} of the weights of its edges.}
of the weights of all its edges. If $G$ has no edges, then $w(G)=1$.
The weight of a set $\SS$ of graphs, $w(\SS)$, is the total weight (the sum of the weights) of its elements; %the graphs belonging to~$\SS$; the weight of the empty set is zero.
$w(\varnothing)=0$. If the weights of all edges are unity, i.e.\ the graphs in $\SS$ are actually unweighted, then $w(\SS)$ reduces to the cardinality of~$\SS$.
The weights of sequences of vertices and edges and of their sets are defined similarly.

For $v_0,v_m\in V(G),$ a $v_0\to v_m$ \emph{path\/} (\emph{simple path}) in $G$ is an alternating sequence of vertices and edges $v_0,\e_1,v_1\cdc\e_m,v_m$ where all vertices are distinct and each $\e_i$ is a $(v_{i-1},v_i)$ edge. The unique $v_0\to v_0$ path is the ``sequence''\,$v_0$ having no edges.
The \emph{length\/} of a path is the number $m$ of its edges. The \emph{weight\/} of a path is the product of the weights of its edges.
The weight of a $v_0\to v_0$ path is~1.

Similarly, a $v_0\to v_m$ \emph{walk\/} (sometimes also called a \emph{route}, cf.\,\cite{Che10bottl}) in $G$ is an \emph{arbitrary\/} alternating sequence of vertices and edges $v_0,\e_1,v_1\cdc\e_m,v_m$ where each $\e_i$ is a $(v_{i-1},v_i)$ edge. %Some (or all) edges in a walk can be loops.
The \emph{length\/} of a walk is the number $m$ of its edges (including loops and repeated edges). The \emph{weight\/} of a walk is the product of the $m$ weights of its edges. The weight of a set of walks is the total weight of its elements. By definition, for any vertex $v_0$, there is one $v_0\to v_0$ walk $v_0$ with length $0$ and weight~1.

We will need several special types of walks. A~\emph{hitting $v_0\to v_m$ walk\/} is a $v_0\to v_m$ walk containing only one occurrence of~$v_m.$
A~$v_0\to v_m$ walk is called a $v_0\to v_0$ \emph{cycle\/} if\xy\footnote{Such a walk is also called a \emph{closed walk}. We use the term \emph{cycle\/} for simplicity; this usage is common in computer science.} $v_m=v_0.$ A~$v_0\to v_0$ cycle is called a \emph{$v_0\tot v_m$ commute cycle\/} if it contains $v_m$ and has no occurrences of $v_0$ strictly between the first appearance of $v_m$ and the final appearance of~$v_0.$

Let $r_{ij}$ be the weight of the set $\RR^{ij}$ of all $i\to j$ walks in $G$, provided that this weight is finite. $R=R(G)=(r_{ij})_{n\times n}$ will be referred to as the \emph{matrix of walk weights}.

By $d^\s(i,j)$ we denote the \emph{shortest path distance}, i.e., the length of a shortest path between $i$ and $j$ in~$G.$
The \emph{weighted shortest path distance\/} $d^\ws(i,j)$ is defined as follows:\footnote{This formula corrects Eq.\,(6.2) in~\cite{KleinRandic93}; cf.\ \cite[Section\:4]{Klein10}.}
\eqss{
\label{e_dws}%+
d^\ws(i,j)=\min_\pi\sum_{\e\xy\in\xy E(\pi)}l_\e,
}
where the minimum is taken over all paths $\pi$ from $i$ to $j$ and the sum is over all edges $\e$ in~$\pi$; $l_\e=1/w_\e$ is sometimes called the \emph{weighted length\/} of the edge $\e$, where $w_\e$ is the weight of this edge (see, e.g.,~\cite{Cinkir11EJC}). In the theory of electrical networks, $l_\e$ is identified with the \emph{resistance\/} of the edge~$\e,$ while $w_\e$ is its conductivity.

\begin{defin}[\!\!\cite{Che10bottl}]
\label{def_trme}
{\rm
Given a graph $G,$ we say that a matrix $S=(s_{ij})\!\in\!\R^{n\times n}$ \emph{determines the transitional measure $s(i,j)\!=\!s_{ij},$ $i,j\!\in\! V,$ for $G$} if $S$
satisfies the \emph{transition inequality\/}\footnote{If $S$ has positive diagonal entries, then the transition inequality is equivalent to $s'_{ij}\,s'_{\!jk}\le s'_{ik},$ where $s'_{ij}
=\frac{s_{ij}}{\sqrt{s_{ii}\xy s_{\xz jj}}},\;$ $i,j,k\in V.$}
\eqs*{
%\label{e_trans}%-
s_{ij}\,s_{\!jk}\le s_{ik}\,s_{\!jj},\quad i,j,k\in V
}
and the \emph{graph bottleneck identity with respect to~$G\!:$}
\eqs*{
%\label{e_boe}%-
s_{ij}\,s_{\!jk}=s_{ik}\,s_{\!jj}
}
holds if and only if all paths in $G$ from $i$ to $k$ contain~$j$.
}
\end{defin}
The transition inequality is a multiplicative analogue of the \emph{triangle inequality for proximities\/} \cite{CheSha97,CheSha98} also called the ``unrooted correlation triangle inequality''~\cite{DezaLaurent97}.

%Suppose that $d(i,j)$ is a distance on~$V(G).$
\begin{defin}[\!\!\cite{Che10bottl}]
\label{d_g-d}
{\rm
For a multigraph $G$ with vertex set $V,$ a function $d\!:V\!\times\!V\to\R$ is called \emph{graph-geodetic\/} provided that
$d(i,j)+d(j,k)=d(i,k)$ holds if and only if every path in $G$ connecting $i$ and~$k$ contains~$j$.
}
\end{defin}

In the following section, we define the class of walk distances and present a number of preliminary results needed in the subsequent study.

\def\baselinestretch{1.0}
\section{The walk distances}
\label{s_prel}

Recall that  $r_{ij}$ is the weight of the set $\RR^{ij}$ of all $i\to j$ walks in $G$ provided that this weight is finite, $R=(r_{ij})_{n\times n}$ being the matrix of walk weights.

For any $t>0,$ consider the graph $G(t)$ obtained from $G$ by multiplying all edge weights by~$t.$
If the matrix $R_t=R(G(t))=(r_{ij}(t))_{n\times n}$ exists,
then\footnote{For an early study of the graph proximity measure $\sum_{k=0}^\infty(tA)^k,$ we refer the reader to \cite{Katz53,Thompson58,Ponstein66,Kasteleyn67,Taylor68}. More recently, it has been explored in \cite{CheSha98,YenSaerensShimbo08,ChelnokovZefirova09,EstradaHigham10}. On counting walks, see also \cite{HararySchwenk79} and on its applications in chemistry,~\cite{KleinPalaciosRandicTrinajstic04}.}
\eq{
\label{e_Rt}%+
R_t=\sum_{k=0}^\infty(tA)^k=(I-tA)^{-1},
}
where $I$ denotes the identity matrix of appropriate dimension.

By assumption, $G$ is connected, while its edge weights are positive, so $R_t$ is positive whenever it exists. Assuming the finiteness of $R_t,$ apply the logarithmic transformation to the entries of~$R_t,$ namely, consider the matrix
\eq{
\label{e_Ha}%+
H_t=\ovec{\ln R_t},
}
where
$\ovec{\varphi(S)}$ stands for elementwise operations, i.e., operations applied to each entry of $S$ separately. Finally, consider the matrix
\eq{
\label{e_Da}%+
D_t=\frac{1}{2}\xy(h_t{\bm1}^\T+\bm1 h_t^\T)-H_t,  %%-H_t^\T),
}
where $h_t$ is the column vector of the diagonal entries of $H_t$ (the \emph{trace vector\/} of~$H_t$), ${\bm1}$ is the vector of ones of appropriate dimension, and $h_t^\T$ and ${\bm1}^\T$ are the transposes of $h_t$ and ${\bm1}$.
An alternative form of \eqref{e_Da} is $D_t=(U_t+U_t^{\xy\T})/2$, where $U_t=h_t{\bm1}^\T-H_t$, and the elementwise form is $d_{ij}(t)=\frac{1}{2}\xy(h_{ii}(t)+h_{\!jj}(t))-h_{ij}(t),\,$ $i,j\in V(G),$ where $H_t=(h_{ij}(t))$ and $D_t=(d_{ij}(t)).$ This is a standard transformation used to obtain a distance from a proximity measure (cf.\ the inverse covariance mapping in~\cite[Section~5.2]{DezaLaurent97} and the cosine law in~\cite{Critchley88}).

In the rest of this section, we present several known facts (lemmas) which will be of use in the study of the walk distances.
The first lemma follows from Theorem~6 in~\cite{Che10bottl}.

\begin{lemma}
\label{l_tm}
For any connected graph $G,$ if the matrix ${R_t=(r_{ij}(t))}$ of walk weights in $G(t)$ exists$,$ then $R_t$ determines a strictly positive transitional measure for~$G.$
\end{lemma}

According to Theorem~1 in~\cite{Che10bottl}, if $S\!=\!(s_{ij})_{n\times n}$ determines a transitional measure for $G$ and has positive off-diagonal entries$,$ then $D=(d_{ij})_{n\times n}$ defined by
$D=\tfrac{1}{2}(h{\bm1}^\T+\bm1 h^\T-H-H^\T),$ where $H=\ovec{\ln S},$ is a matrix of distances on\/~$V(G).$ Moreover, by Theorem~2 in~\cite{Che10bottl} this distance is graph-geodetic. Along with Lemma~\ref{l_tm} this implies the following lemma, which appears in ~\cite{Che10bottl} as item~2 of Corollary~2.

\begin{lemma}
\label{l_di}
For any connected $G,$ if ${R_t=(r_{ij}(t))}$ exists$,$ then the matrix $D_t=(d_{ij}(t))$ defined by~\eqref{e_Rt}--\eqref{e_Da} determines a graph-geodetic distance $d_t(i,j)=d_{ij}(t)$ on~$V(G).$
\end{lemma}

\begin{defin}
\label{d_walkD}
{\rm
For a connected graph $G,$ the \emph{walk distances\/} on~$V(G)$ are the functions $d_t(i,j)\!:V(G)\!\times\! V(G)\to\R$ and the functions positively proportional to them$,$ where $d_t(i,j)=d_{ij}(t)$ and $D_t=(d_{ij}(t))$ is defined by~\eqref{e_Rt}--\eqref{e_Da}.
}
\end{defin}

Regarding the finiteness of $R_t,$ since for a connected graph, $A$ is irreducible, the Perron-Frobenius theory of nonnegative matrices provides the following result (cf.~\cite[Theorem~4]{Thompson58}).
\begin{lemma}
\label{l_finite}
For any weighted adjacency matrix $A$ of a connected graph $G,$ the series
$R_t=\sum_{k=0}^\infty(tA)^k$ with $t>0$ converges to $(I-tA)^{-1}$ if and only if\/ $t<\rho^{-1},$ where $\rho=\rho(A)$ is the spectral radius of\/~$A.$ Moreover$,$ $\rho$ is an eigenvalue of~$A;$ as such $\rho$ has multiplicity~$1$ and a positive eigenvector.
\end{lemma}

Eigenvalue $\rho=\rho(A)\xy$ is called the \emph{Perron root\/} of~$A.$
If $x$ is an eigenvector of $A$ associated with~$\rho,$ then the probability vector $p=x/\|x\|_1$ is called the \emph{Perron vector of}~$A.$

\begin{lemma}
\label{l_expD}
For any vertices $\,i,j\in V(G)$ and $\,0<t<\rho^{-1},$
\eqss{
\label{e_d1}%+
d_t(i,j)
=-\ln\biggl(\xz\frac{r_{ij}(t)}{\sqrt{r_{ii}(t)\,r_{\!jj}(t)}}\xz\biggr).
}
\end{lemma}

Lemma~\ref{l_expD} is a corollary of \eqref{e_Ha} and \eqref{e_Da} (cf.\ Eq.~(11) in~\cite{Che10bottl}).
The author is grateful to Michel Deza for mentioning the genetic distance by Nei~\cite{Nei72}, which has a form similar to~\eqref{e_d1}.

Lemma~\ref{l_decoR} appeared in~\cite{Che10bottl} as Eq.\,(23). Despite its simplicity, it plays an important role in the subsequent study.
\begin{lemma}
\label{l_decoR}
If the matrix $R=(r_{ij})$ exists$,$ then for any vertices $i,j\in V(G),$
\eq{
\label{e_rij}%+
r_{ij}=r_{ij(1)}\,r_{\!jj},
}
where $r_{ij(1)}=w(\RR^{ij(1)})$ is the weight of the set $\,\RR^{ij(1)}$ of all $i\to j$ hitting walks in~$G$.
\end{lemma}

\section{Two expressions for the walk distances}
\label{s_trans}

The first result enables one to interpret the walk distances in terms of specific walks in~$G$.
Technically, it is a consequence of the previous lemmas.
\begin{thm}
\label{t_commu}
For any $t\in\,]0,\rho^{-1}[,$ the matrix of walk distances $D_t$ has the representation
\eq{
\label{e_disCyc}%+
D_t
=-\frac{1}{2}\,\ovec{\ln(R_{(1)t}R_{(1)t}^\T)}
=-\frac{1}{2}\,\ovec{\ln C^{\ms\circlearrowright}_t},
}
where $R_{(1)t}=(r_{ij(1)}(t))_{n\times n}$ is the matrix of hitting walk weights in~$G(t),$ $C^\circlearrowright_t=(c_{ij}^\circlearrowright(t))_{n\times n},$ and $c_{ij}^\circlearrowright(t)=w(\CC^{i\tot j}_t)$ is the weight of the set $\,\CC^{i\tot j}_t$ of all\, $i\tot j$ commute cycles in~$G(t)$.
\end{thm}

\proof
By Lemma~\ref{l_finite}, if $0<t<\rho^{-1},$ then the distance matrix $D_t$ exists and by Lemmas~\ref{l_expD} and~\ref{l_decoR},
for any vertices $i$ and $j$ we have
\eqss{
\label{e_d2}%+
d_t(i,j)
=-\frac{1}{2}\ln\!\frac{r^2_{ij}(t)}
                       {r_  {ii}(t)\,
                        r_{\!jj}(t)}
=-\frac{1}{2}\ln\xz\left(r_{  ij(1)}(t)\,
                         r_{\!ji(1)}(t)\xz\right).
}
Observing now that there is a natural bijection between $\RR^{ij(1)}_t\times\RR^{ji(1)}_t$ and $\CC^{i\tot j}_t$ we obtain %that
\eq{
\label{e_r-c}%+
c_{ij}^\circlearrowright(t)
=r_{ij(1)}(t)\,r_{\!ji(1)}(t),\quad i,j\in V(G),
}
which implies~\eqref{e_disCyc}.
\qed\bigskip

By virtue of Theorem~\ref{t_commu}, there is a certain analogy between the walk distances and the classical commute time distance.
One of the consequences of Theorem~\ref{t_commu} is that $w(\CC^{i\tot j}_t)<1$  whenever $R_t$ exists and $i\ne j.$
A ``topological'' interpretation of the walk distances is presented in~\cite{CheDeza11}.

The following result provides an expression for the walk distances which will be of use in the sequel.

\begin{thm}
\label{t_decom}
For any connected $G,$ any vertices $i,j\in V(G),$ and any $t\in\,]0,\rho^{-1}[,$
\eqs*{
%\label{e_disDeco}%-
d_t(i,j)
=-\frac{1}{2}\ln\bigl(\xz
 (t^{-1}\xz I-A_\jj)^{-1}_ia\_{\j j}\,
 (t^{-1}\xz I-A_\ii)^{-1}_ja\_{\i i}
  \bigr),
}
where $M_\jj$ is the submatrix of $M$ obtained by the removal of row $j$ and column $j,$ $M^{-1}_i$ is the $i$th row of $M^{-1},$ and $a\_{\j j}$ is the $j$th column of~$A$ with $a_{jj}$ removed.
\end{thm}

\proof
Theorem~\ref{t_decom} is immediate from~\eqref{e_d2} and the following lemma.
\qed
\begin{lemma}
\label{l_rij1}
In the notation of Theorems~$\ref{t_commu}$ and\/ $\ref{t_decom},$
\eq{
\label{e_rij1}%+
r_{ij(1)}(t)=(t^{-1}\xz I-A_\jj)^{-1}_ia\_{\j j}
}
whenever $t\in\,]0,(\rho(A_\jj))^{-1}[.$
\end{lemma}

\proof
Observe that any $i\to j$ hitting walk in $G(t)$ can be uniquely decomposed into: (1)~some $i\to k$ walk in the subgraph $G_{\j\xy}(t)$ of $G(t)$ obtained by the removal of vertex $j$ and all edges incident to it and (2)~a $(k,j)$ edge. If $0<t<(\rho(A_\jj))^{-1},$ then the total weights of the $i\to k$ walks in $G_{\j\xy}(t)$ form the $i$th row of $(I-tA_\jj)^{-1}$, whereas the total weights of the $(k,j)$ edges (with $k\ne j$) form the vector~$ta\_{\j j}.$ The desired expression follows.
\qed\medskip

When considering graph distances, of major interest are the proportions of distances for different pairs of vertices rather than the distances themselves. On the other hand, for studying the limit properties, it is convenient to consider, among the positive multiples of $d_t(i,j)$ (see Definition~\ref{d_walkD}), the specific walk distances $d^\W_\a(i,j)$ with: ``{\footnotesize W}'' referring to ``walk'',
\eq{
\label{e_gam1}%+
d^\W_\a(i,j)=\th\,d_t(i,j),
}
$\a$ being the parameter connected with both $t$ and $\rho$ by
\eq{
\label{e_param}%+
\a=(t^{-1}-\rho)^{-1},
}
and $\th$\/ being the scaling factor given by
\eq{
\label{e_gam2}%+
\th=\ln(\ee+\a^{\frac{2}{n}})\frac{\a-1}{\ln\a}.
}
The factor $\th$ as a function of $\a$ and~$n$ is assumed to extend to $\a=1$ by continuity:
\eqs*{
%\label{e_gam3}%-
\th=\ln(\ee+1)\;\;\text{whenever}\;\;\a=1.
}
These parameterization and scaling will prove convenient in the following sections. In particular, it is worth mentioning that they ensure comparability of the walk distances with the logarithmic forest distances \cite{Che08dist} (cf.\ Section~\ref{s_forest}).

\section{The short walk distance}
\label{s_short}

Consider the behavior of the walk distances $d_\a^\W(i,j)$ as $\a\to 0^+$ ($t\to 0^+$).
The corresponding limit of $d_\a^\W(i,j)$ (if it exists and provides a distance) can be termed the \emph{short walk distance\/} because $t\to 0^+$ leads to neglecting long walks in \eqref{e_Rt}. It turns out that the short walk distance coincides with the classical shortest path distance~$d^\s(i,j).$

\begin{thm}
\label{t_short}
For any vertices $i,j\in V,$
\eqs*{
%\label{e_limt0}%-
\lim_{\a\to0^+}d_\a^\W(i,j)=d^\s(i,j),
}
where $d^\s(i,j)$ is the shortest path distance between $i$ and $j$ in~$G.$
\end{thm}

\proof
For any vertices $i$ and $j\!\ne\! i,$ let $m\!=\!d^\s(i,j).$ Let $r_{ij}^{(m)}$ be the $ij$-entry of~$A^m.$ Using Lemma~\ref{l_expD} and \eqref{e_gam1}--\eqref{e_gam2} yields
\eqs*{
%\label{e_pr1}%-
   \lim_{\a\to0^+}d_\a^\W(i,j)
&=&\lim_{\a\to0^+}\frac{1}{2}\ln(\ee+\a^{\frac{2}{n}})\frac{\a-1}{\ln\a}\ln\!\frac{r_{  ii}((\rho+\a^{-1})^{-1})\,
                                                                                   r_{\!jj}((\rho+\a^{-1})^{-1})}
                                                                                {r^2_{  ij}((\rho+\a^{-1})^{-1})}\nonumber\\
&=&-\frac{1}{2}\lim_{\a\to0^+}(\ln\a)^{-1}\ln\!\frac{r_{ii}(\a)\,r_{\!jj}(\a)}  {r^2_{  ij}(\a)}\nonumber\\
&=&-\frac{1}{2}\lim_{\a\to0^+}(\ln\a)^{-1}\ln\!\frac{(1+o(1))(1+o(1))}{(\a^mr^{(m)}_{ij}+o(\a^m))^2}\nonumber\\
&=&            \lim_{\a\to0^+}(\ln\a)^{-1}(m\ln\a +\ln r^{(m)}_{ij})
 = m,
 }
where $o(f(\a))$ are such terms that $\frac{o(f(\a))}{f(\a)}\to0$ as $\a\to0^+.$ This completes the proof.
\qed\bigskip

\section{The long walk distance}
\label{s_long}

Consider the asymptotic behavior of the walk distances as $\a\to\infty$ ($t\to(\rho^{-1})^-$). First, the behavior of $R_t=\sum_{k=0}^\infty(tA)^k$ is clear from the following lemma.

\begin{lemma}
\label{l_resolv}
For any connected graph $G,$
\eq{
\label{e_limpo}%+
 \lim_{t\to(\rho^{-1})^-}(t^{-1}-\rho)R_t
=\lim_{\a\to\infty}\a^{-1}R_t
=\rho\xy\tilde p\xy\tilde p^{\xy\T},
}
where $\tilde p=p/\|p\|_2,$ $\,p=(p_1\cdc p_n)^\T$ is the Perron vector of~$A,$ and \eqref{e_param} is used.
\end{lemma}

\proof
Eq.\,\eqref{e_limpo} can be easily derived from, say, Theorem~3.1 in~\cite{Meyer74} and the fact that $\tilde p\xy\tilde p^{\xy\T}$ is the eigenprojection of $A$ corresponding to~$\rho$ (see also~\cite{Rothblum76,Rothblum81}). When applying this theorem, to verify that the limit in \eqref{e_limpo} exists, one should observe that the index of $A$ at $\rho$ is~$1$ since $A$ is diagonalizable as a Hermitian matrix.
\qed\medskip

While the entries of $R_t$ (the total weights of walks between vertices) tend to infinity as $t\to(\rho^{-1})^-,$ the weights of hitting walks and commute cycles remain finite.

\begin{corol}[of Lemma~\ref{l_resolv}]
In the notation of Theorem~$\ref{t_commu},$ for any vertices $i,j\in V,$
\label{c_rinf}
\eqss{
\label{e_rinf}%+
\lim_{t\to(\rho^{-1})^-}r_{ij(1)}(t)
&=&\frac{p_i}{p_j},\\
\label{e_c-1}%+
\lim_{t\to(\rho^{-1})^-}c_{ij}^\circlearrowright(t)
&=&1.
}
\end{corol}

\proof
Combining \eqref{e_rij} and \eqref{e_limpo} yields
\eqs*{
 \lim_{t\to(\rho^{-1})^-}r_{ij(1)}(t)
=\lim_{t\to(\rho^{-1})^-}\xz\frac{r_{ij}(t)}{r_{\!jj}(t)}
=\frac{p_ip_j}{p_j^2}=\frac{p_i}{p_j}.
}
In view of \eqref{e_r-c}, Eq.\:\eqref{e_c-1} also holds.
\qed

\medskip
It follows from~\eqref{e_disCyc} and~\eqref{e_c-1} that the distances $d_{ij}(t)$ vanish as $t\to(\rho^{-1})^-$ in spite of the infiniteness of~$R_t.$ Furthermore, since $A$ is irreducible, $\rho(A_\jj)<\rho(A)$ for any $j\in V$ %\cite[item\:5.5.6 of Part\:II]{MarcusMinc64}.
\cite[Ch.\:III, \S\:3.4]{Gantmacher59-2}.
Consequently, by \eqref{e_rij1} and \eqref{e_r-c}, $\rho I -A_\jj$ is non-singular, $r_{ij(1)}(\rho^{-1})$ and $c_{ij}^\circlearrowright(\rho^{-1})$ make sense, and so Eqs.\,\eqref{e_rinf} and~\eqref{e_c-1} should be supplemented by
\eqss{
\label{e_rinf1}%+
r_{ij(1)}(\rho^{-1})
&=&\frac{p_i}{p_j},\\
\label{e_c-11}%+
c_{ij}^\circlearrowright(\rho^{-1})
&=&1.
}

Substituting \eqref{e_c-11} in~\eqref{e_disCyc} results in $D_{\rho^{-1}}=0.$ However, by~\eqref{e_param}, $\a$ is indefinite at $t=\rho^{-1},$ so $d^\W_\a(i,j)$ defined by~\eqref{e_gam1} is indefinite as well. Therefore $\lim_{t\to(\rho^{-1})^-}d_\a^\W(i,j)=\lim_{\a\to\infty}d_\a^\W(i,j)$ is worth evaluating. Let us study this limit.

\medskip
%A~vis-a-vis of the short walk distance is
We define the \emph{long walk distance\/} $\dLW(i,j)$ as follows:
\eq{
\label{e_dLW}%+
\dLW(i,j)
=\lim_{\a\to\infty}d_\a^\W(i,j),\quad i,j\in V,
}
provided that the limit exists and induces a distance function.
In Theorem~\ref{t_LWd}, we obtain a closed formula for $\dLW(i,j)$; after that we give it an interpretation and discuss two examples.

\begin{thm}
\label{t_LWd}
For any vertices $i,j\in V$ such that $i\ne j,$
\eq{
\label{e_LWd}%+
\dLW(i,j)
=n^{-1}\Bigl(p_i^{-1}(\L_\jj)^{-1}_i\xy p\_\j
            +p_j^{-1}(\L_\ii)^{-1}_j\xy p\_\i\Bigr),
}
where $\L=\rho I-A,\xy$ $p=(p_1\cdc p_n)^\T$ is the Perron vector of $A,$ $ p\_\j$ is $p$ with\/ $p\_j$ removed\/$,$ and the other notation is the same as in Theorem~$\ref{t_decom}.$
\end{thm}

\proof
Using Theorem~\ref{t_decom}, \eqref{e_gam1}--\eqref{e_gam2} and the Taylor expansion we obtain
\eqs*{
\dLW(i,j)
&=& \lim_{\a\to\infty}d_\a^\W(i,j)\nonumber\\
&=&-\lim_{\a\to\infty}\frac{1}{2}\ln(\ee+\a^\frac{2}{n})\frac{\a-1}{\ln\a}\ln\bigl((t^{-1}\xz I-A_\jj)^{-1}_ia\_{\j j}\,
                                                                                   (t^{-1}\xz I-A_\ii)^{-1}_ja\_{\i i}\bigr)\nonumber\\
&=&-n^{-1}\!\lim_{\a\to\infty}\a\ln      \bigl([(\L_\jj)^{-1}
                                        -\a^{-1}(\L_\jj)^{-2}+o(\a^{-1})]_i\,a\_{\j j}\pha\nonumber\\
&&                        \hspace{5.85em}\times[(\L_\ii)^{-1}
                                        -\a^{-1}(\L_\ii)^{-2}+o(\a^{-1})]_j\,a\_{\i i}\bigr)\nonumber\\
&=&-n^{-1}\!\lim_{\a\to\infty}\a\ln  \bigl(p_{ij}p_{ji}
                              -\a^{-1}[\xy Y\!(j)_   i Y\!(j)a\_{\j j}
                                           p_{ji}
                                          +p_{ij}\xy
                                           Y\!(i)_{\!j}Y\!(i)a\_{\i i}\xy]+o(\a^{-1})\bigr),\pha\nonumber
}
where
\eq{
\label{e_Yi}%+
Y\!(i)=(\L_\ii)^{-1},\quad i\in V
}
and
\eq{
\label{e_pij}%+
p_{ij}=Y\!(j)_i\,a\_{\j j},\quad i,j\in V,\;\; i\ne j.
}
To proceed, we need the following lemma.
\begin{lemma}
\label{l_pij}
\eq{
\label{e_pij1}%+
p_{ij}
=\frac{p_i}{p_j},\quad i,j\in V(G),\;\; i\ne j,
}
where $p_{ij}$ is defined by~\eqref{e_pij} and $p=(p_1\cdc p_n)^\T$ is the Perron vector of~$A.$
\end{lemma}

Lemma~\ref{l_pij} can be proved using \eqref{e_rij1} and \eqref{e_rinf1}, however, it is more instructive to give a direct proof.

\medskip
\proof
As $p$ is the Perron vector of $A,$ it obeys $\L p=0.$ Removing the $j$th equation from this linear system and rearranging all the $p_j$'s to the right side yields $\L_\jj\xy p\_\j=p_ja\_{\j j}.$
Therefore, since $\L_\jj$ is non-singular, $p_j^{-1} p\_\j=Y\!(j)\xy a\_{\j j}$ holds, as required.
\qed

\medskip
Using Lemma~\ref{l_pij} we now complete the proof of Theorem~\ref{t_LWd}:
\eqs*{
\dLW(i,j)
&=&-n^{-1}\!   \lim_{\a\to\infty}\ln\bigl(1-\a^{-1}\!\bigl[\,p_i^{-1}{p\_j}\,Y\!(j)\_   i \,p_j^{-1} {p\_\j}
                                                            +p_j^{-1}{p\_i}\,Y\!(i)\_{\!j}\,p_i^{-1} p\_\i\,\bigr]\!\bigr)^\a\pha\nonumber\\
&=&-n^{-1\ms}                    \ln\exp\bigl(-p_i^{-1}\xy Y\!(j)_   i \,p\_\j
                                              -p_j^{-1}\xy Y\!(i)_{\!j}\,p\_\i\bigr)_{\ms}\pha\nonumber\\
&=& n^{-1\ms}            \Bigl(p_i^{-1}(\L_\jj)^{-1}_   i p\_\j
                              +p_j^{-1}(\L_\ii)^{-1}_{\!j}p\_\i\Bigr)^\mz,\nonumber
}
as desired.
\qed\bigskip

The symmetric irreducible
singular M-matrix $\L=\rho I-A$ plays a central role in this paper. It can be termed the \emph{para-Laplacian matrix\/} of~$G.$
$\L$~has rank $n-1$ and is positive semidefinite. %HornJohnson "Topics..." p. 131, problem 27.

Expression~\eqref{e_LWd} can be written in a more elegant form.
\begin{corol}[of Theorem~\ref{t_LWd}]
\label{c_LWd0}
For any vertices $i,j\in V$ such that $i\ne j,$
\eq{
\label{e_LWd0}%+
\dLW(i,j)
=n^{-1}\Bigl[(\rho I-B_\jj)^{-1}_i
            +(\rho I-B_\ii)^{-1}_j\Bigr]\bm1,
}
where $B=P^{-1}AP$ and\/ $P=\diag p.$
\end{corol}

The proof of Corollary~\ref{c_LWd0} is straightforward.
It should be noted that $Q=(\rho P)^{-1}AP$ is a stochastic matrix which can be naturally attached to~$G.$ In terms of $Q,$ one can write:
\eqs*{
%\label{e_dLWQ}%-
\dLW(i,j)
=(n\rho)^{-1}\bigl[(I-Q_\jj)^{-1}_i
                  +(I-Q_\ii)^{-1}_j\bigr]\bm1,\quad j\ne i.
}

Now let us give an interpretation of $\dLW(i,j)$ in terms of walks.

Denote by ${\cal C}^{i(j)}$ the set of all cycles (i.e., closed walks) in~$G(\rho^{-1})$ that
\begin{itemize}
\item start and finish at $i;$
\item consist of two consecutive walks such that the first one does not contain $j$ and finishes at some vertex $k$ which is marked; the second one does not contain $i,$ except for its end vertex.
\end{itemize}
Let $c^{i(j)}=w({\cal C}^{i(j)})$ be the weight of the set ${\cal C}^{i(j)}.$

\begin{corol}[of Theorem~\ref{t_LWd}]
\label{c_LWd}
For any vertices $i,j\in V$ such that $i\ne j,$
\eqs*{
%\label{e_LWd1}%-
\dLW(i,j)
=(n\rho)^{-1}(c^{i(j)}+c^{j(i)}).
}
\end{corol}

\proof
Let $y_{ik}^{(j)}$ be the $k$th element of $Y\!(j)_i=(\L_\jj)^{-1}_i,$ where $k\in V\!\smallsetminus\xz\{j\}.$
It follows from the proof of Lemma~\ref{l_rij1} (Section~\ref{s_trans}) that $\rho y_{ik}^{(j)}$ is equal to the total weight of the $i\to k$ walks in $G_{\j\xy}(\rho^{-1}).$
Using Theorem~\ref{t_LWd} and \eqref{e_rinf1} we obtain
\eqss{
\dLW(i,j)
 &=&n      ^{-1}\Bigl(\sum_{k\in V,\,k\ne j}     y_{ik}^{(j)}p_i^{-1}p_k
                     +\sum_{k\in V,\,k\ne i}     y_{jk}^{(i)}p_j^{-1}p_k\Bigr)\nonumber\\
 &=&(n\rho)^{-1}\Bigl(\sum_{k\in V,\,k\ne j}\rho y_{ik}^{(j)}\,r_{ki(1)}\xz(\rho^{-1})
                     +\sum_{k\in V,\,k\ne i}\rho y_{jk}^{(i)}\,r_{kj(1)}\xz(\rho^{-1})\Bigr)\nonumber\\
 &=&(n\rho)^{-1}(c^{i(j)}+c^{j(i)}),
\label{e_terpr1}%+
}
as required.
\qed

\medskip
By virtue of Corollary~\ref{c_LWd}, it can be said that $\dLW(i,j)$ is proportional to the sum of
the weights of certain walks starting at $i,$ avoiding $j$, and then returning to $i$ and
               certain walks starting at $j,$ avoiding $i$, and then returning to~$j.$

\begin{corol}[of Theorem~\ref{t_LWd}]
\label{c_LWd2}
For any connected graph $G,$ the function $\dLW(i,j)$ is a metric.
\end{corol}

\proof
Since $\dLW(i,j)$ is a finite limit of distances, it suffices to prove that ${\dLW(i,j)\ne0}$ when\-ever $j\ne i.$ This follows from the non-emptiness of ${\cal C}^{i(j)}$ for all $i$ and $j\ne i$ and Corollary~\ref{c_LWd}.
\qed

\begin{example}
{\rm For the unweighted path $P_4$ (Fig.\;\ref{f_1}), we find using Theorem\;\ref{t_LWd} that $\dLW(1,2)/\dLW(2,3)=(1+\sqrt{5})/2,$ the \emph{golden ratio.}

\medskip
\begin{figure}[ht!]
\begin{center}
{\small 1 --------- 2 --------- 3 --------- 4\hspace{4em} 1 --------- 2 --------- 3 --------- 4 --------- 5}
\end{center}

\hspace{9.5em}$P_4$\hspace{17em}$P_5$
\caption{The paths $P_4$ and $P_5$.\label{f_1}}
\end{figure}

In general, it can be shown that the long walk distance between central adjacent vertices in a path is smaller than that between peripheral adjacent vertices. For example, for $P_5$ (Fig.~\ref{f_1}), ${\dLW(1,2)/\dLW(2,3)=2.}$ This distinguishes the long walk distance (and all walk distances) from the logarithmic forest distances (cf.\ the remark {``On the `mixture' of the shortest-path and resistance distances''} in Section~6 of~\cite{Che08dist}). Since the long walk distance is the limit of graph-geodetic distances, the distances between non-adjacent vertices $i$ and $j$ in a path are equal to the sum of distances between the subsequent vertices in the subpath connecting $i$ and~$j.$
}
\end{example}

\begin{example}
{\rm Consider the weighted paths with four and five vertices and the edge weights shown in Fig.~\ref{f_2}.

\medskip
\begin{figure}[ht!]
\noindent
{\scriptsize$\;\;$
{\footnotesize Weights:}\hspace{1.8em}$\sqrt{2}$\hspace{4.5em}1\hspace{4.5em}$\sqrt{2}$\hspace{4.9em}
{\footnotesize Weights:}\hspace{1.8em}$\sqrt{2}$\hspace{4.6em}1\hspace{4.7em}1\hspace{4.4em}$\sqrt{2}$}

\vspace{-1.3em}
\begin{center}
{\small 1 --------- 2 --------- 3 --------- 4\hspace{5em} 1 --------- 2 --------- 3 --------- 4 --------- 5}
\end{center}

\hspace{9.2em}$P_4$\hspace{17.7em}$P_5$
\caption{Two weighted paths, $P_4$ and $P_5$.\label{f_2}}
\end{figure}

The results are as follows: for $P_4,\,$ $\dLW(1,2)=\dLW(2,3)=0.75;\,$ for $P_5,\,$ $\dLW(1,2)=\dLW(2,3)=0.8.$ The same pattern is preserved for all weighted paths of this kind. Say, for $P_{10}$ with the two terminal weights $\sqrt{2}$ and the other weights~$1,$ all the long walk distances between adjacent vertices are~$0.9$ (and $\frac{n-1}n$ for $P_n,$ $n\!>\!2$). Thus, the weights of $\sqrt{2}$ completely \emph{compensate\/} the ``extremality'' of path's terminal vertices with respect to the long walk distance.
}
\end{example}

\def\baselinestretch{0.99}
\section{The $\bm\ee$-walk distances which generalize\\ the weighted shortest path distance}
\label{s_e-modi}

In Section~\ref{s_prel}, the graph $G(t)$ was constructed by multiplying all edge weights in $G$ by~$t.$
Now consider a more sophisticated transformation:\xy\footnote{Obviously, for~\eqref{e_trex}, $\lim_{\a\to0^+}w(\a)=0$ and  $\lim_{\a\to\infty}w(\a)=w/\rho$ as well as for the transformation ${w(\a)=tw=w/(\rho+\a^{-1})}$ which we used earlier (see~\eqref{e_param}). However, the \emph{rate\/} of convergence as a function of $w$ for these two transformations cardinally differs.}
\eq{
\label{e_trex}
w(\a)
=\tfrac{w}{\rho}\xy\ee^{-\frac{1}{\a w}},\quad \a>0,
}
where $w$ is any edge weight in $G,$ $w(\a)$ is the weight of the corresponding edge in the transformed graph $\ti G(\a),$ and $\rho$ is the Perron root of the weighted adjacency matrix $A$ of~$G.$
The total edge weights $a_{ij}(\a)=\sum_{q=1}^{n_{ij}}w_{ij}^q(\a)$ for all pairs of vertices form the weighted adjacency matrix $A(\a)$ of the transformed graph $\ti G(\a)$ (cf.\:\eqref{e_aij}).
The matrix of walk weights of this graph, $\ti R_\a=(\ti r_{ij}(\a))_{n\times n},$ provided that it exists, has the representation
\eq{
\label{e_Rt1}%+
\ti R_\a=\sum_{k=0}^\infty A(\a)^k=(I-A(\a))^{-1}.
}
When the series in \eqref{e_Rt1} converges, we define the \emph{modified\/} walk distance $\ti d^\W_\a(i,j)$ by means of
\eq{
\label{e_Ha1}%+
\ti H_\a=\th_{\xz\a}\xy\a\,\ovec{\ln\ti R_\a},
}
where $\th_{\xz\a}$ is a positive scaling factor,
\eq{
\label{e_Da1}%+
\ti D_\a=\frac{1}{2}\xy(\ti h_\a{\bm1}^\T+\bm1 \ti h_\a^\T)-\ti H_\a, %-\ti H_\a^\T),
}
where $\ti h_\a$ is the column vector of the diagonal entries of $\ti H_\a$ (cf.\ \eqref{e_Ha} and~\eqref{e_Da}), and
\eq{
\label{a_dWmod}%+
\ti d^\W_\a(i,j)=\ti d_{ij}(\a),
}
where $(\ti d_{ij}(\a))_{n\times n}=\ti D_\a.$

$\ti d^\W_\a(i,j)$ are graph-geodetic distances on $V(G)$ since Lemmas~\ref{l_tm} and~\ref{l_di} remain valid for~$\ti R_\a$.

\begin{defin}
{\rm
For a connected $G,$ the \emph{$\ee$-walk distances\/} on~$V(G)$ are the functions $\deW(i,j)\stackrel{\rm def}{=}\ti d^\W_\a(i,j)$ defined by \eqref{e_trex}--\eqref{a_dWmod}. More generally, a \emph{modified walk distance\/} is a distance that fits within the framework of \eqref{e_Rt1}--\eqref{a_dWmod}
with some edge weight transformation~$w(\a).$
}
\end{defin}

The following expression for $\deW(i,j)$ is analogous to the representation~\eqref{e_d1} of the walk distances. It is easily obtained by combining \eqref{e_trex}--\eqref{e_Da1}.

\begin{lemma}
\label{l_modif1}
For any $\a>0,$
\eqs*{
%\label{e_dwm1}%-
\deW(i,j)
=-\th_{\xz\a}\xy\a\ln\biggl(\xz\frac{\ti r_{ij}(\a)}{\sqrt{\ti r_{ii}(\a)\,\ti r_{\!jj}(\a)}}\xz\biggr),\quad i,j\in V,\;i\ne j,
}
where
\eqs*{
%\label{e_dwm2}%-
\ti r_{ij}(\a)
=\sum_{\om\in\RR^{ij}}\rho^{-m_\om}w\_\om\ee^{-d_\om/\a},
\quad i,j\in V,
}
$m\_\om$ and $w\_\om$ are the length and the weight of the walk\/ $\om,$ respectively$,$ $\,d_\om=\sum_{\e\in E(\om)}l_\e,$ $\,l_\e=w_\e^{-1},$ and $E(\om)$ is the multiset of the edges of\/ $\om$ $(\om$ may have repeated edges$).$
\end{lemma}

\medskip
For the
$\ee$-walk distances, an analogue of Theorem~\ref{t_decom} holds (and has a similar proof).
\begin{lemma}
\label{l_modilo}
For a connected $G,$ $\,i,j\in V(G),$ and any $\a>0,$ in the notation of Theorem~$\ref{t_decom},$
\eq{
\label{e_modilo}%+
\deW(i,j)
=-\frac{\th_{\xz\a}\xy\a}{2}\ln\bigl(\xz
 (I-A(\a)_\jj)^{-1}_ia(\a)\_{\j j}\,
 (I-A(\a)_\ii)^{-1}_ja(\a)\_{\i i}
  \bigr).
}
\end{lemma}

Lemmas~\ref{l_modif1} and~\ref{l_modilo} are used in the proof of the following theorem describing the limiting properties of the $\ee$-walk distances (which differ from those of $d_\a^\W(i,j)$).
Suppose that for the $\ee$-walk distances $\deW(i,j),$ $\th_{\xz\a}$ is such that
\eqss{
\label{e_th_inf}
\lim_{\a\to0^+}\th_{\xz\a}=1\;\;\;\text{and}\;\;\,\lim_{\a\to\infty}\th_{\xz\a}=\th_\infty\in\R_+.
}

\begin{thm}
\label{t_modif1}
For any vertices $i,j\in V$ such that $j\ne i,$
\eqs*{
%\label{e_lima0}%-
\lim_{\a\to0^+}\deW(i,j)=d^\ws(i,j),
}
where $d^\ws(\cdot,\cdot)$ is the weighted shortest path distance\/ \eqref{e_dws} and
\eq{
\label{e_lima8}%+
\lim_{\a\to\infty}\deW(i,j)
=\frac{\th_\infty}{2}\Bigl(p_i^{-1}\bigl(\xz\xz(\L_\jj)^{-1}\cA_\j\bigr)_{\xz i}
                          +p_j^{-1}\bigl(\xz\xz(\L_\ii)^{-1}\cA_\i\bigr)_{\!  j}\Bigr)p,
}
where $\L=\rho I-A,$ $p$ is the Perron vector of $A,$ $\cA=(\Check a_{ij})_{n\times n}$ results from $A$ by replacing every nonzero entry by\/ $1,$ and $\cA_\i$ is $\cA$ with the $i$th row removed.
\end{thm}

\proof
Using Lemma~\ref{l_modif1} and \eqref{e_th_inf}, for any vertices $i$ and $j\ne i$ we obtain
\eqss{
\label{e_lima0p1}%+
     \lim_{\a\to0^+}\deW(i,j)
 =-\!\lim_{\a\to0^+}\!\a\ln\!\frac{\ti r_{ij}(\a)}{\sqrt{(1+o(1))(1+o(1))}}
 =-\!\lim_{\a\to0^+}\!\a\ln       {\ti r_{ij}(\a)}.
}

Observe that if $\om,\om'\in\RR^{ij}$ and, in the notation of Lemma~\ref{l_modif1}, $d_{\om'}<d_\om,$ then for all sufficiently small $\a>0$,
$\rho^{-m_{\om'}}w\_{\om'}\ee^{-d_{\om'}/\a}
>\rho^{-m_ \om  }w\_ \om  \ee^{-d_ \om  /\a}$ holds. Consequently, there exists $\a_0>0$ such that for all $\a\in\,]0,\a_0[$
and some $\ka_{ij}(\a)$ satisfying $1\le\ka_{ij}(\a)\le|\RR^{ij}|$,

\eq{
\label{e_kappa2}%+
\ti r_{ij}(\a)
=\ka_{ij}(\a)\,\rho^{-m_\omb}w\_\omb\ee^{-d_\omb/\a}
}
is true, where $\omb\in\RR^{ij}$ is a walk such that either              \newline
(a) $d_{\omb}<d_{\om}$ or                                                \newline
(b) $d_{\omb}=d_{\om}$ and $\rho^{-m_\omb}w\_\omb\ge\rho^{-m_\om}w\_\om$ \newline
holds w.r.t.\ all $\om\in\RR^{ij}$. By definition \eqref{e_dws}, in this case, $d_{\omb}=d^\ws(i,j)$. Using \eqref{e_lima0p1} and \eqref{e_kappa2} we obtain
$$
  \lim_{\a\to0^+}\deW(i,j)
=-\lim_{\a\to 0^+}\a\big(\ln(\ka_{ij}(\a)\,\rho^{-m_\omb}w\_\omb)-d_{\omb}/\a\big)
 =d^\ws(i,j).%\eqno{\qed}
$$

Now we prove~\eqref{e_lima8}. Using a different parameterization of the function \eqref{e_trex}:
\eq{
\label{e_bex}%+
w(\a)=\ti w(\g)
=\tfrac{w}{\rho}\xy\ee^{-\frac{\g}{w}},
}
where $\g=\a^{-1},$ observe that
\eq{
\label{e_bexd}%+
\ti w'(0)=-\rho^{-1},
}
where $\ti w'(\g)$ is the derivative of $\ti w(\g)$ with respect to~$\g.$
Denote by $\ti A(\g)$ the weighted adjacency matrix of the graph modified through~\eqref{e_bex}.
As $\a\to\infty$ ($\g\to0^+$), Eqs.\;\eqref{e_bex} and \eqref{e_bexd} and the definition of $\cA$ yield
\eqss{
\label{e_pertA}%+
A(\a)_\jj
=      \ti A(0)_\jj-\a^{-1}\rho^{-1}\cA_\jj +o(\a^{-1})
=\rho^{-1}(A   _\jj-\a^{-1}\!       \cA_\jj)+o(\a^{-1}),\quad j\in V.
}
For the vector $a(\a)\_{\j j}$ (the $j$th column of $A(\a)$ with $a(\a)_{jj}$ removed) this implies that
\eqss{
\label{e_perta}%+
a(\a)\_{\j j}
=\rho^{-1}(a\_{\j j}-\a^{-1}\Check a\_{\j j}+o(\a^{-1})),\quad j\in V.
}

Substituting \eqref{e_pertA} and \eqref{e_perta} in \eqref{e_modilo} and denoting by $\Check a\_{\i i}$ the $i$th column of $\cA$ with
$\Check a\_{ii}$ removed result in
\eqs*{
%\label{e_lima8p}%-
              \lim_{\a\to\infty}\deW(i,j)
&=&-\th_\infty\lim_{\a\to\infty}\frac{\a}{2}\ln
\bigl(
                      [I-\rho^{-1}(A_\jj-\g\cA_\jj)+o(\g)]^{-1}_i\rho^{-1}(a\_{\j j}-\g\Check a\_{\j j}+o(\g))\\
&&\hspace{5.8em}\times[I-\rho^{-1}(A_\ii-\g\cA_\ii)+o(\g)]^{-1}_j\rho^{-1}(a\_{\i i}-\g\Check a\_{\i i}+o(\g))
\bigr)\\
&=&-\th_\infty\lim_{\a\to\infty}\frac{\a}{2}\ln
\bigl(
                      [\L_\jj+\g\cA_\jj+o(\g)]^{-1}_i(a\_{\j j}-\g\Check a\_{\j j}+o(\g))\\
&&\hspace{5.8em}\times[\L_\ii+\g\cA_\ii+o(\g)]^{-1}_j(a\_{\i i}-\g\Check a\_{\i i}+o(\g))
\bigr).
}
Observe that when $\g\to0^+,$
\eqs*{
(\xz(\L_\ii+\g\cA_\ii)-\g\cA_\ii)
    (\L_\ii+\g\cA_\ii)^{-1}
=I-\g\cA_\ii
    (\L_\ii)^{-1}+o(\g)
}
is true, from which
\eqss{
\label{e_pert}%+
 (\L_\ii+\g       \cA_\ii)^{-1}
= Y\!(i)-\g Y\!(i)\cA_\ii\xy Y\!(i)+o(\g),\quad j\in V
}
holds, where $Y\!(i)=(\L_\ii)^{-1}$ (see \eqref{e_Yi}).

Using \eqref{e_pert}, \eqref{e_pij} and \eqref{e_pij1} and denoting by $p\_\i$ the Perron vector $p$ of $A$ with $p\_i$ removed, we can now complete the proof:
\eqs*{
              \lim_{\a\to\infty}\deW(i,j)
&=&-\th_\infty\lim_{\a\to\infty}\frac{\a}{2}\ln
\bigl(
                      [\xy Y\!(j)-\g Y\!(j)\cA_\jj\xy Y\!(j)]_i\xy(a\_{\j j}-\g\Check a\_{\j j})\\
&&\hspace{5.9em}\times[\xy Y\!(i)-\g Y\!(i)\cA_\ii\xy Y\!(i)]_j\xy(a\_{\i i}-\g\Check a\_{\i i})+o(\g)\!
\bigr)\\
&=&-\th_\infty\lim_{\a\to\infty}\frac{1}{2}\ln
\Bigl(\frac{p_i}{p_j}\frac{p_j}{p_i}
        -\a^{-1}\Bigl([\xy Y\!(j)\cA_\jj\xy]_i\, p_j^{-1}p\_\j\xy\frac{p_j}{p_i}+Y\!(j)_     i \xy\Check a\_{\j j}\xy\frac{p_j}{p_i}\\
&&\hspace{11.4em}    +[\xy Y\!(i)\cA_\ii\xy]_j\, p_i^{-1}p\_\i\xy\frac{p_i}{p_j}+Y\!(i)_{\xz j}\xy\Check a\_{\i i}\xy\frac{p_i}{p_j}
          \Bigr)+o(\a^{-1})\!
\Bigr)^\a\\
&=&\frac{\th_\infty}{2}
\bigl(p_i^{-1}[\xy Y\!(j)\cA_\j\xy]_i
     +p_j^{-1}[\xy Y\!(i)\cA_\i\xy]_j\bigr)p,
}
which coincides with the desired expression.
\qed

\medskip
Is there any connection between the limiting $\ee$-walk distance \eqref{e_lima8} and the long walk distance $\dLW(i,j)$ defined by~\eqref{e_dLW}? Let
\eqss{
\label{e_MLW}
d^\LeW(i,j)=
\lim_{\a\to\infty}\deW(i,j),\quad i,j\in V(G)
}
(L{\large $\ee$}{\xz}W is the abbreviation for ``long $\ee$-walk'').
In fact, $d^\LeW(i,j)$ is a distance, which is guaranteed by Theorem~\ref{t_LW2}. Prior to formulating this theorem, we provide a ``topological'' interpretation of $d^\LeW(i,j).$

\smallskip
Recall the interpretation of the long walk distance $\dLW(i,j)$ given by Corollary~\ref{c_LWd} (Section\:\ref{s_long}) in terms of specific cycles in~$G(\rho^{-1}).$ Such a cycle belonging to ${\cal C}^{i(j)}$ is an $i$-to-$i$ cycle that consists of two consecutive walks: the first one does not contain $j$ and finishes at some vertex $k$ which is marked; the second walk does not contain $i,$ except for its end vertex. Let us take such a cycle and remove the edge connecting $k$ with the subsequent vertex in the cycle. Let $\ti{\cal C}^{i(j)}$ be the set of resulting sequences and let $\ti c^{\,i(j)} =w(\ti{\cal C}^{i(j)})$ be its weight. Each element of $\ti{\cal C}^{i(j)}$ can be treated as a ``cycle with a jump''. Indeed, one can imagine a point moving along the cycle, reaching $k,$ and then jumping to the next vertex instead of traversing the edge leading to it.

\begin{corol}[of Theorem~\ref{t_modif1}]
\label{c_LeWd}
For any vertices $i,j\in V$ such that $i\ne j,$
\eqs*{
%\label{e_LeWd1}%-
d^\LeW(i,j)
=\frac{\th_\infty}{2\rho}\,(\ti c^{\,i(j)}+\ti c^{\xy j(i)}).
}
\end{corol}

Similarly to the long walk distance, the long $\ee$-walk function
$d^\LeW(i,j)$ is large when the set comprising specific $i\to i$ cycles avoiding, on the first stage,~$j$ along with
                                               specific $j\to j$ cycles avoiding, on the first stage,~$i$ is ``heavy''.

\medskip
\proof
Using Theorem~\ref{t_modif1} we obtain (cf.\ the proof of Corollary~\ref{c_LWd}):
\eqss{
d^\LeW(i,j)
&=&\frac{\th_\infty}{2\rho}\Bigl(\sum_{k,\xy q\in V,\,k\ne j}\rho y_{ik}^{(j)}\,\Check a_{kq}\,r_{qi(1)}\xz(\rho^{-1})
                                +\sum_{k,\xy q\in V,\,k\ne i}\rho y_{jk}^{(i)}\,\Check a_{kq}\,r_{qj(1)}\xz(\rho^{-1})\Bigr)\nonumber\\
&=&\frac{\th_\infty}{2\rho}\,(\ti c^{\,i(j)}+\ti c^{\xy j(i)}),
\label{e_terpr2}%+
}\nopagebreak
as required.
\qed

\begin{thm}
\label{t_LW2}
In the notation of\/ $\eqref{e_dLW},$ $\eqref{e_th_inf}$ and $\eqref{e_MLW},$ and Theorem\/~$\ref{t_modif1},$ if
\eqss{
\label{e_thin}%+
\th_\infty=
\frac{2}{n}\cdot\frac{p^\T\xz(A/\rho)p}{p^\T\!\cA p},
}
then for all vertices $i,j\in V,$
\eqs*{
%\label{e_lima8a}%-
d^\LeW(i,j)=\dLW(i,j).
}
\end{thm}

\begin{remark}
{\rm
Observe that $\frac{p^\T\xz(A_{\phantom\ms}/\rho)p}{p^{\T^{\phantom i}}\!\!\!\cA p}$ is the weighted average, with weights $p_ip_j,$ of the nonzero entries $a_{ij}/\rho$ of $A/\rho.$
Since, by assumption\:\eqref{e_th_inf}, $\lim_{\a\to0^+}\th_\a=1,$ a scaling factor $\th_\a$ in~\eqref{e_Ha1} that ensures $d^\LeW(i,j)=\dLW(i,j)$ can be defined, for instance, as follows:
\eqs*{
%\label{e_tha}%-
\th_\a
=\frac{\th_\infty\xy\a+\beta}{\a+\beta},
}
where $\th_\infty$ is given by~\eqref{e_thin} and $\beta$ is a positive parameter.
}
\end{remark}

{\noindent\bf Proof of Theorem~\ref{t_LW2}.}
Let
\eqs*{
\eta_{kq}^{ij}
=\d_{kj}\xy\rho\xy y_{ik}^{(j)}\,r_{qi(1)}\xz(\rho^{-1})
+\d_{ki}\xy\rho\xy y_{jk}^{(i)}\,r_{qj(1)}\xz(\rho^{-1}),\quad i,j,k,q\in V,
}
where
$y_{ik}^{(j)}$ is the $k$th element of $Y\!(j)_i=(\L_\jj)^{-1}_i$ and
\eqs*{
\d_{kj}=
\begin{cases}
1,&k\ne j,\\
0,&k  = j.
\end{cases}
}
It follows from~\eqref{e_terpr1} and~\eqref{e_terpr2} that for any vertices $i$ and $j\ne i,$
\eqss{
\label{e_divdd}%+
\frac{\dLW(i,j)}{d^\LeW(i,j)}
=\frac{2}{n\xy\th_\infty}\cdot
\frac{\suml_{k,\xy q\in V}\eta_{kq}^{ij}\,       a_{kq}\xy\rho^{-1}}
     {\suml_{k,\xy q\in V}\eta_{kq}^{ij}\,\Check a_{kq}}.
}

Using \eqref{e_rinf1} one can represent the vector $\eta_{\ast q}^{ij}=(\eta_{1q}^{ij}\cdc\eta_{nq}^{ij})^\T$ in the form
\eqss{
\label{e_eta}%+
\eta_{\ast q}^{ij}
=\rho\Bigl(
 \oY(j)^i\xy\frac{p_q}{p_i}
+\oY(i)^j\xy\frac{p_q}{p_j}
\Bigr)
=\frac{\rho\xy p_q}{p_ip_j}\bigl(
 \oY(j)^i\xy p_j
+\oY(i)^j\xy p_i\bigr),
\quad i,j,q\in V,
}
where $\oY(j)$ is $Y\xz(j)$ supplemented by row $j$ and column $j$ consisting of zero entries and $\oY(j)^i$ is the $i$th column of $\oY(j).$
The rest of the proof is based on the following lemma.

\smallskip
\begin{lemma}
\label{l_aver}
Vector $\oY(j)^i\xy p_j
       +\oY(i)^j\xy p_i$ with $j\ne i$ is a positive multiple of the Perron vector~$p.$
\end{lemma}

\smallskip
\proof
Performing multiplication of block matrices and using \eqref{e_Yi}--\eqref{e_pij1} one can verify that   $\L\oY(j)^i\xy p_j$ is the vector whose $i$th element is $p_j,$ the $j$th element is $-p_i,$ and the other elements are zero. Similarly, in the vector $\L\oY(i)^j\xy p_i,$ the $j$th element is $p_i,$ the $i$th element is $-p_j,$ and the remaining elements are zero.
Therefore, ${\L\bigl(\oY(j)^i\xy p_j
                    +\oY(i)^j\xy p_i\bigr)=0}$ and so
                    $\oY(j)^i\xy p_j
                    +\oY(i)^j\xy p_i$ is a positive (as $\dLW(i,j)>0$) multiple of vector $p$ spanning $\Ker\L.$
\qed

\medskip
By \eqref{e_eta} and Lemma~\ref{l_aver}, every vector $\eta_{\ast q}^{ij}$ is proportional to~$p.$ Owing to the factor $p_q$ in~\eqref{e_eta}, every row of the matrix $(\eta_{kq}^{ij})_{n\times n}$ indexed by $k$  and~$q$ is proportional to~$p^\T.$ Therefore
\eq{
\label{e_matnu}%+
(\eta_{kq}^{ij})_{n\times n}
=\mu_{ij}\xy p\xy p^\T,
}
where $\mu_{ij}$ is a factor of proportionality. Substituting \eqref{e_matnu} in \eqref{e_divdd} leads to the result.
\qed

\def\baselinestretch{1.0}
\section{Logarithmic forest distances as a subclass\\ of walk distances}
\label{s_forest}

For a graph $G$ and a parametric family of functions $\varphi_\a\!:\R_+\to\R_+,$ $\,\a\in\A\subseteq\R,$ consider the matrices
\eq{
\label{QaO}%+
Q_\a=(I+L_\a)^{-1},
}
where $L_\a\!=\!\diag(A_\a\bm1)\!-\!A_\a$ and $A_\a$ are the Laplacian and weighted adjacency matrices of the graph $G_\a$ that differs from $G$ by the edge weights only:
\eqs*{
%\label{e_phi}%-
w_\a=\varphi_\a(w)
}
for any edge weight $w$ in $G$ and the corresponding weight $w_\a$ in~$G_\a.$

The \emph{logarithmic forest distances\/} on $V(G)$ determined by the parametric edge weight transformation $\varphi_\a$ are obtained~\cite{Che08dist} from the matrices $Q_\a$ through the familiar conversions
\eq{
\label{HaO}%+
H_\a=\th\,\ovec{\ln Q_\a},
}
where $\th$ is a positive scaling factor generally depending on $\a$ and~$G$ and
\eq{
\label{DaO}%+
D_\a=\frac{1}{2}(h_\a{\bm1}'+\bm1 h'_\a)-H_\a.
}

The simplest edge weight transformation $\varphi_\a(w)=\a w,$ $\a>0$ determines~\cite{Che08dist} a specific family of logarithmic forest distances whose limiting cases are the shortest path distance and the resistance distance.

In this section, we establish a connection between the walk distances and the logarithmic forest distances.

Let us say that \emph{$\ti G$ is a balance-graph of $G$} if $\ti G$ is obtained from $G$ by attaching some loops and assigning the loop weights that provide $\ti G$ with uniform weighted vertex degrees. More formally, ${V(\ti G)\!=\!V(G),}$ ${E(G)\!\subseteq\! E(\ti G),}$ the edges in $E(G)$ have the same weights in~$\ti G,$ $\,{E(\ti G)\xz\smallsetminus\xz E(G)}$ is comprised of loops, and $A(\ti G)$ has constant row sums.
\emph{Balancing $G$ by loops\/} will mean constructing any balance-graph of~$G$.

\begin{thm}
\label{t_forest}
For any connected graph $G,$ the family of logarithmic forest distances \eqref{QaO}--\eqref{DaO} with any edge weight transformation $\varphi_\a(w)$ coincides with a certain family of modified walk distances \eqref{e_Rt1}--\eqref{a_dWmod} obtained through balancing the graphs $G_\a$ by loops.
\end{thm}

\proof
For each $\a\in\A,$ choose any
\eq{
\label{e_mal}%+
m_\a\ge\max\xy\bigl(\l_{ii}(\a)\mid i\in V(G)\xz\bigr),
}
where $\l_{ii}(\a)$ are the diagonal entries of~$L_\a=L(G_\a).$ Since $G$\/ is connected, $m_\a>0.$ Set
\eqss{
\label{e1_LtoA3}%+
A(\a)
=(m_\a+1)^{-1}(m_\a I-L_\a).
}

Obviously, $A(\a)$ defined by \eqref{e1_LtoA3} is the weighted adjacency matrix of the graph with loops $\ti G_\a$ obtained from $G$\/ by transforming the edge weights in accordance with
\eqs*{
%\label{e1_LtoAp}%-
\ti w(\a)
=(m_\a+1)^{-1}\varphi_\a(w),
}
attaching a loop to each vertex $i$ such that $m_\a>\l_{ii}(\a),$ and assigning the loop weights ${(m_\a+1)^{-1}(m_\a-\l_{ii}(\a))}$; such weights provide $A(\a)$ with constant row sums $\,m_\a/(m_\a+1).$
Thereby, the $\ti G_\a$'s are obtained from the $G_\a$'s through \emph{balancing by loops}.

The Perron root of $A(\a),$ $\,m_\a/(m_\a+1),$ is less than~1. Consequently,
\eqss{
\label{e1_1}%+
\ti R_\a\stackrel{\rm def}{=}
\sum_{k=0}^\infty A(\a)^k
=\bigl(I-A(\a)\bigr)^{-1}
}
is a finite matrix of walk weights in~$\ti G_\a.$

Substituting \eqref{e1_LtoA3} into \eqref{e1_1} yields
\eqss{
\label{e1_LtoA1}%+
\ti R_\a
=(m_\a+1)(I+L_\a)^{-1}
=(m_\a+1)Q_\a.
}

Passing $Q_\a$ through the conversions \eqref{HaO}--\eqref{DaO} leads to the logarithmic forest distance with parameter~$\a.$
Passing $\ti R_\a$ through the same conversions fits within the framework of \eqref{e_Rt1}--\eqref{e_Da1} and so it generates a modified walk distance with parameter~$\a.$

Finally, observe that the multiplier $(m_\a+1)$ in \eqref{e1_LtoA1} does not survive the conversions \eqref{HaO}--\eqref{DaO}, so the two above distances coincide. %Moreover, these two distances are generated with the same parameter~$\a.$
Considering the whole domain $\A$ of $\a$ leads us to recognize that the initial family of logarithmic forest distances and the family of modified walk distances constructed by means of $\ti R_\a$ and \eqref{HaO}--\eqref{DaO} also coincide.
\qed\medskip

The following corollary applies to the simplest case of logarithmic forest distances in which ${\varphi_\a(w)=\a w}$ and so $L_\a=\a L$ \cite[Section\:2]{Che08dist}, where $L=(\l_{ij})=L_1$ ($\a=1$) is the Laplacian matrix of~$G.$
\begin{corol}
\label{c_ForSimp}
For any connected graph $G,$ if\/ $\varphi_\a(w)=\a w,$ then the family of logarithmic forest distances \eqref{QaO}--\eqref{DaO} with $\A=\R_+$ and\/ $\th$ given by~\eqref{e_gam2} coincides with the family of walk distances~\eqref{e_gam1} calculated for any balance-graph of\/~$G$.
\end{corol}

\proof
Following the lines of the proof of Theorem~\ref{t_forest}, denote $m_1$ by $m$ (see\:\eqref{e_mal}).
Since $w_\a=\a w,$ for every $\a>0$ we have $\max\xy\bigl(\l_{ii}(\a)\mid i\!\in\! V(G)\xz\bigr)=\a\max\xy\bigl(\l_{ii}(1)\mid i\!\in\! V(G)\xz\bigr)\le\a m.$ So in\:\eqref{e_mal} one can set
\eqss{
\label{e2_0}%+
m_\a=\a m,\quad\a>0.
}

Define $A(\a)$  by \eqref{e1_LtoA3} and \eqref{e2_0} and let
\eqss{
\label{e2_1}%+
\ti A_\a
=\a^{-1}(\a m+1)\xy A(\a).
}
Then
\eqss{
\label{e2_3}%+
\ti A_\a=mI- L,\quad
\rho(\ti A_\a)=m,
}
where $\rho(\ti A_\a)$ is the Perron root of~$\ti A_\a.$ Thus, $\ti A_\a$ does not depend on $\a$; denote it by~$\ti A.$

Substituting \eqref{e2_1} into \eqref{e1_1} yields
\eqs*{
%\label{e2_2}%-
\ti R_\a
=\bigl(I-\a\xy(\a m+1)^{-1}           \ti A\xy\bigr)^{-1}
=\bigl(I-\bigl(\rho+\a^{-1}\bigr)^{-1}\ti A\xy\bigr)^{-1},
}
where $\rho=\rho(\ti A).$
Setting, in accordance with \eqref{e_param}, $\,t=(\rho+\a^{-1})^{-1}$ we have
\eqs*{
%\label{e2_7}%-
\ti R_\a
=\bigl(I-t\ti A\xy\bigr)^{-1},
}
which coincides with \eqref{e_Rt} for the graph $\ti G$ whose weighted adjacency matrix is~$\ti A.$
Passing $\ti R_\a$ through the conversions \eqref{HaO}--\eqref{DaO} with $\th$ given by \eqref{e_gam2} provides exactly the walk distance~\eqref{e_gam1} with parameter~$\a.$ Passing $Q_\a$ through the same conversions results in the logarithmic forest distance under consideration, which, by Theorem~\ref{t_forest}, coincides with the above walk distance. Since this holds for every $\a\in\A,$ the two families of distances coincide.

Finally, observe that by \eqref{e2_3},
\eq{
\label{e_Atil}
\ti A=mI-\diag(A\bm1)+A,
}
thus, $\ti G$ can be constructed by attaching a loop to each vertex $i$ such that $m>\l_{ii}$ and assigning the loop weights that provide $\ti G$ with uniform weighted vertex degrees. Obviously, each balance-graph of $G$ can be obtained in this way. The corollary is proved.
\qed

\section{Connections between long walk distance\\ and resistance distance}
\label{s_walkresi}
\subsection{Resistance distance as the long walk distance\\ in a balance-graph}

It follows from Corollary~\ref{c_ForSimp} that the logarithmic forest distances in $G$ with edge weight transformation $\varphi_\a(w)=\a w$ coincide with the walk distances in any balance-graph of~$G.$ Since by Proposition\:3 in~\cite{Che08dist}, the resistance distance is a limiting case of the logarithmic forest distances, the resistance distance can be obtained within the framework of walk distances.

\begin{corol}[of Theorem\:\ref{t_forest}]
\label{c_Resist}
For any connected $G,$ the resistance distance in $G$ coincides with the long walk distance $\dLW(i,j)$ defined by~\eqref{e_dLW} in $\ti G,$ where $\ti G$ is any balance-graph of\/~$G.$
\end{corol}

Corollary\;\ref{c_Resist} is immediate from Proposition\;3 in \cite{Che08dist} and Corollary\;\ref{c_ForSimp}. It enables one to apply to the resistance distance any result obtained for the long walk distance. In particular, Corollary\;\ref{c_LWd} of Section\;\ref{s_long} (with $\rho=m,$ where $m$ is the uniform weighted vertex degree of a balance-graph of~$G$) provides a kind of topological interpretation of the resistance distance, whereas Theorem~\ref{t_LWd} gives the following expression.

\begin{corol}
\label{c_Resist1}
For any connected graph $G$ on $n$ vertices$,$ let $L$ be the Laplacian matrix of~$G$ and let $\dr(\cdot,\cdot)$ be the resistance distance on~$V(G).$ Then for any $i,j\in V(G)$ such that $j\ne i,$
\eqs*{
\dr(i,j)
=n^{-1}\bigl(\xz(L_\jj)^{-1}_i
               +(L_\ii)^{-1}_j\bigr)\bm1
}
holds$,$ where $\bm1$ is the vector of\/ $n-1$ ones and\/
$(L_\jj)^{-1}_i$ is the $i$th row of the inverse principal submatrix~$L_\jj.$
\end{corol}

\proof
By Corollary~\ref{c_Resist}, the resistance distance in $G$ coincides with the long walk distance in any balance-graph $\ti G$ of~$G.$ The weighted adjacency matrix $\ti A$ (see \eqref{e_Atil}) of the balance-graph $\ti G$ having all weighted vertex degrees $m$ is a nonnegative irreducible matrix with row sums~$m.$ Therefore $\rho(\ti A)=m$ and $p(\ti A)=n^{-1}\bm1,$ where $\rho(\ti A)$ and $p(\ti A)$ are the Perron root and Perron vector of $\ti A,$ respectively. Substituting this and \eqref{e_Atil} into the expression for $\dLW(i,j)$ given by Theorem~\ref{t_LWd} yields the desired equation.
\qed\medskip

Note that Corollary~\ref{c_Resist1} can also be proved using the results of~\cite{BapatSivasubramanian11}.
\bigskip

It follows from the proof of Corollary\:\ref{c_ForSimp} that the logarithmic forest distance with parameter $\a$ coincides with the walk distance \eqref{e_gam1}, provided that $\a$ is defined by \eqref{e_param} and the graph has been balanced by loops. This justifies the reparameterization~\eqref{e_param}.

Attaching the ``balancing loops'' leads to a model with a uniform connection resource possessed by all vertices: a lack of external connections is filled up by self-connections. As has been seen in this section, in such models, the logarithmic forest distances appear. These treat two peripheral adjacent vertices in a path as being closer to each other \cite{Che08dist} than two central adjacent vertices are. It was mentioned in the Introduction that friendship is one of the relationships for which such a model can be considered. It may be appropriate when several people have a similar combined resource of friendship + self-absorption, but they are not equal in their ability to make friends. In contrast to this, the examples of Sections\;\ref{s_long} and\;\ref{s_example} demonstrate that the walk distances are able to treat central adjacent vertices in a path as being closer to each other than the peripheral adjacent vertices are, which also may be relevant to certain applications.

\subsection{Long walk distance as the resistance distance\\
            in a modified graph}

The connection between long walk distance and resistance distance is two-way. Namely, the following relationship supplements Corollary\:\ref{c_Resist}.

\begin{thm}
\label{t_rel}
Let $G$ be a connected graph on $n$ vertices with weighted adjacency matrix $A.$ Suppose that $p$ is the Perron vector of~$A,\xy$ $p'=\sqrt{n}p/\|p\|_2,$ and $P'=\diag p'.$ Then the long walk distance in $G$ coincides with the resistance distance in the graph $G'$ whose weighted adjacency matrix is~$P'AP'.$
\end{thm}

We first prove two lemmas.

\begin{lemma}
\label{l_newRep}
Let $X=\diag x,$ where $x\in\R^n$ is non-negative and $x\ne\bm0.$ Then in the notation of Theorem\/\:$\ref{t_rel}$ and Lemma\/\:$\ref{l_aver},$
\eqs*{
%\label{e_LWd0}%-
\dLW(i,j)
=\frac{\|p\|_2^2}{n\xy p^\T x}
\Bigl[(P^{-1}\oY(j)\xy X)_i
    + (P^{-1}\oY(i)\xy X)_j\Bigr]\bm1,\quad j\ne i
}
holds$,$ where $P=\diag p.$
\end{lemma}

\proof
For $j\ne i,$ we have
\eqs*{
\frac{\|p\|_2^2}{n\xy p^\T x}
\Bigl[(P^{-1}\oY(j)\xy X)_i
    + (P^{-1}\oY(i)\xy X)_j\Bigr]\bm1
&=&\frac{\|p\|_2^2}{n\xy p^\T x}   \bigl(\xy p_i^{-1}\oY(j)_i+
                                             p_j^{-1}\oY(i)_j\bigr)x\nonumber\\
&=&\frac{\|p\|_2^2}{n\xy p^\T x}(p_ip_j)^{-1}\xz\bigl(p_j\oY(j)^i+
                                                      p_i\oY(i)^j\bigr)^{\xz\T}x.
}
By Lemma\:\ref{l_aver}, $p_j\oY(j)^i+p_i\oY(i)^j=\beta_{ij}\xy p$ holds$,$ where $\beta_{ij}>0$ is a factor of proportionality. Consequently, Theorem\:\ref{t_LWd} yields
\eqs*{
\frac{\|p\|_2^2}{n\xy p^\T x}
\Bigl[(P^{-1}\oY(j)\xy X)_i
    + (P^{-1}\oY(i)\xy X)_j\Bigr]\bm1
&=&\frac{\|p\|_2^2}{n\xy p^\T x}(p_ip_j)^{-1}\frac{\beta_{ij}\xy p^\T x}
                                                  {\beta_{ij}\xy p^\T p}\bigl(p_j\oY(j)_i+
                                                                              p_i\oY(i)_j\bigr)p\\
&=&\frac{\|p\|_2^2}{n\xy p^\T x}\!\cdot\!\frac{p^\T x}{\|p\|_2^2}\bigl(  p_i^{-1} Y\xz(j)_i\,p\_\j+
                                                                         p_j^{-1} Y\xz(i)_j\,p\_\i\bigr)\\
&=&\dLW(i,j),
}
and the lemma is proved.
\qed

\begin{remark}
\label{r_Ytil}
{\rm
It can be shown that in Lemma\:\ref{l_newRep}, $\oY(j)$ can be replaced by $\tY(j)=(\rho I-A_{\jjt\xz})^{-1},\,$ where $A_{\jjt\xz}$
is $A$ with the entries in the $j$th row and $j$th column replaced by zero.}
\end{remark}

Lemma\:\ref{l_newRep} provides one more formula for the long walk distance, which can be computationally cheaper than~\eqref{e_LWd0}.

\begin{corol}
\label{c_simpLW}
In the notation of Lemma\/\:$\ref{l_newRep}$ and Remark\/\:$\ref{r_Ytil},$
\eqs*{
%\label{e_LWd0}%-
\dLW(i,j)
=\frac{\|p\|_2^2}{n}
\Bigl[((\rho I-A_{\jjt\xz})P)^{-1}_i
    + ((\rho I-A_{\iit\xz})P)^{-1}_j\Bigr]\bm1,\quad j\ne i.
}
\end{corol}

\proof
This follows from Lemma\:\ref{l_newRep} by setting $x=\bm1,$ which implies $X=I.$
\qed

\begin{lemma}
\label{l_BapLap}
Let $A(G')=\beta PA(G)P,$ where $G$ is connected$,$ $P=\diag p,$ $p$ is the Perron vector of~$A(G),$ and\/ $\beta>0.$ Then $L(G')=\beta P\L(G)P,$ where $\L(G)=\rho I-A(G)$ and $\rho$ is the Perron root of~$A(G).$
\end{lemma}

\proof
Obviously, the non-diagonal entries of $L(G')$ coincide with those of $\beta P\L(G)P.$ Finally, $\beta P\L(G)P$ has zero row sums$:$
$\beta P\L(G)P\bm1=\beta P(\rho I-A(G))p=\beta P(\rho p-\rho p)=\bm0.$
\qed
\medskip\medskip

{\noindent\bf Proof of Theorem\:\ref{t_rel}.}
Setting $x=(p_1^{-1}\cdc p_n^{-1})^\T$ and using Corollary\:\ref{c_Resist1} and Lemmas\:\ref{l_BapLap} and\:\ref{l_newRep} we have
\Up{.38}
\eqss{
d_{G'}^{\xy\om}(i,j)
&=&n^{-1}                          \bigl(   \bigl(L'_\jj                    \bigr)^{-1}_i+
                                            \bigl(L'_\ii                    \bigr)^{-1}_j\bigr)\bm1\nonumber\\
&=&\frac{\|p\|_2^2}{n^2}           \bigl[\xz\bigl((P     \L P )_\jj  \bigr)^{-1}_i+
                                            \bigl((P     \L P )_\ii  \bigr)^{-1}_j\bigr]\bm1\nonumber\\
&=&\frac{\|p\|_2^2}{n\xy p^\T x}   \bigl[\xz\bigl( P^{-1}\oY(j)\xy X\bigr)_i+
                                            \bigl( P^{-1}\oY(i)\xy X\bigr)_j\bigr]\bm1\nonumber\\
&=&d_G^{\xy\LW}\xz(i,j),
}
as desired.
\qed

\medskip
Theorem\:\ref{t_rel} enables one to utilize all facts and expressions known for the resistance distance to calculate and study the long walk distance. In particular, Corollaries\;\ref{c_geoEuc} and\;\ref{c_frResi} follow.

\begin{corol}%[of Theorem\:\ref{t_rel}]
\label{c_geoEuc}
The long walk distance is graph-geodetic\/$;$ it is a squared Euclidean distance.
\end{corol}

\proof
Note that the resistance distance has these properties \cite{Gurvich10,MooreSubak-Sharpe68} and use Theorem\:\ref{t_rel}.
\qed

\begin{corol}%[of Theorem\:\ref{t_rel}]
\label{c_frResi}
In the notation of Theorem\/\:$\ref{t_rel},$
\eqss{
\label{e_LWd01}%+
\dLW(i,j)
&=&(-1)^{u+v}\frac{\det\xy\xy(L'_\ii)_\jj}{\det L'_{\uv}},\quad j\ne i,\quad\forall\, u,v\in V,\\
\label{e_LWd02}%+
\dLW(i,j)
&=&\l'^{\bm-}_{ii}+\l'^{\bm-}_{jj}-2\l'^{\bm-}_{ij},\\
\label{e_LWd03}%+
\dLW(i,j)
&=&x_\u^\T(i,j)(L'_\vu)^{-1}x\_\v(i,j),\quad j\ne i,\quad\forall\, u,v\in V.
}
where $L'=P'\L\xy P',\,$ $L'^{\bm-}=(\l'^{\bm-}_{ij})$ is any g-inverse\/\footnote{$Z$ is a \emph{g-inverse\/} \cite{RaoMitra71} of $X$ whenever $X=XZX.$} of\:$L',$ and\/ $x(i,j)$ is the $n$-vector whose $i$th element is $+1,$ $j$th element is\, $-1,$ and the other elements are~$0$.
\end{corol}

\proof
This follows from Theorem\:\ref{t_rel}, Lemma\:\ref{l_BapLap}, and three classical expressions for the resistance distance (see \cite[Eq.\,(17)]{SharpeSpain60}, \cite[Theorem~7--4]{SeshuReed61}, and \cite[Eq.\,(14)]{SharpeStyan65} %and \cite[Eq.\,(13)]{SharpeStyan67}
for \eqref{e_LWd01}, {\cite[Eq.\,(13)]{SharpeStyan67}} and \cite[Theorem\:10.1.4]{RaoMitra71} for \eqref{e_LWd02}, and \cite[Eqs.\,(14)--(15)]{SharpeSpain60}, \cite[the first part of Eq.\,(16)]{SharpeStyan65}, and \cite[Eq.\,(15)]{SharpeStyan67} for \eqref{e_LWd03}; cf.\ \cite[Chapter\:7]{SeshuReed61} and \cite{KleinRandic93,StyanSubak-Sharpe97,Bapat99RD}).
\qed
\medskip

For any symmetric irreducible Laplacian matrix $L,$ a simple choice of $L^{\bm-}$ is ${(L+\bar J)^{-1},}$ where $\bar J=\frac1n\bm{11}^\T$ \cite[Section\:10.1.3]{RaoMitra71}. Another choice is the group inverse ${L^\#=(L+\bar J)^{-1}-\bar J}$ (which for $L$ is also the Moore-Penrose generalized inverse~$L^+$). The latter formula, due to Sharpe and Styan \cite{SharpeStyan65b} (see also \cite[Theorem\:10.1.2]{RaoMitra71} and \cite[Propositions\:15,\:16]{CheAga02ap}), has been rediscovered several times.
Alternatively, $\l^\#_{ij}=\frac1{n^2}\bm1^\T L_\ji^{-1}\bm1,$ ${i,j=\1n}$~\cite{SharpeStyan65b}. The general form of $L^{\bm-}$ is $L^\#+a\bm1^\T+\bm1b^{\hy\T},$ where $a$ and $b$ are arbitrary $n$-vectors~\cite{SharpeStyan65}. In particular, $-\frac12D$ is a g-inverse of $L,$ where $D$ is the matrix of resistance distances corresponding to $L$ \cite{SharpeStyan65,Youla59}.

\medskip
Finally, let us mention three simple expressions for $\dLW(i,j)$ in terms of $\L$ obtained in~\cite{CheBapBal11}.

\begin{thm}[\!\!\cite{CheBapBal11}]
\label{t_LW_I-A}
In the notation of Theorem\/\:$\ref{t_rel},$ for all $i,j\in V$ such that $j\ne i,$
\eqs*{
%\label{e_dLW_L1}%-
\dLW(i,j)
&=&\frac{\det\xy\xy(\L_{\xy\ii})_\jj}{p'^2_j\det \L_{\xy\ii_{\ms}}},\\
%\label{e_dLW_L2}%-
\dLW(i,j)
&=&z^\T(i,j)\xy\L^{\bm-}z(i,j)^\mz,\\
%\label{e_dLW_L3}%-
\dLW(i,j)
&=&z_{\xy\u}^\T\xy(i,j)\hy(\L_\vu)^{-1}z\_\v(i,j)^\mz,\quad \forall\, u,v\in V,
}
where
$\L=\rho I-A,\xy$ $\L^{\bm-}$ is any g-inverse of~$\L,$ and\/ $z(i,j)$ is the $n$-vector whose $i$th element is $1/p'_i,$ $j$th element is\, $-1/p'_j,$ and the other elements are\/~$0$.
\end{thm}

%\begin{remark}
%{\rm
Since for any balance-graph $\ti G$ of $G,\,$ $p'(\ti G)\!=\!\bm1,\,$ $\L(\ti G)\!=\!L(G),$ while $L(G)$ is an equicofactor matrix, Theorem\:\ref{t_LW_I-A} generalizes the three classical expressions for the resistance distance reproduced in Corollary\:\ref{c_frResi}.

More generally, the long walk distance can be considered as the counterpart of the resistance distance obtained by replacing the Laplacian matrix $L\!=\!\diag(A\bm1)\!-\!A$ and the vector~$\bm1$ which spans $\Ker L$ with the ``\emph{para-Laplacian}'' matrix $\L=\rho I-A$ and the vector $p'$ spanning $\Ker\L.$ If $G$ is balanced, i.e., $A$ has constant row sums, then these distances coincide.
%}
%\end{remark}

\section{Several metrics on the path of length~$\bm3$}
\label{s_example}

\begin{figure}[ht!]
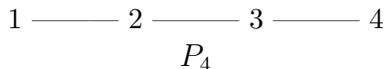

\begin{center}
{\small 1 --------- 2 --------- 3 --------- 4}

$P_4$
\end{center}
\Up{.7}
\caption{The path $P_4.$\label{f_3}}
\end{figure}

The simplest graph on which the difference between the new and classical metrics can be illustrated is the path on $4$ vertices (Fig.\;\ref{f_3}).

\bigskip

Some properties of different metrics on $P_4$ are summarized in Table~\ref{ta_di}.

\begin{table}[ht]
$$
\begin{tabular}{|l||>{$}r<{$}|>{$}r<{$}|>{$}r<{$}|}                                                                             \hline
Metric, $d$                          &\frac{d\ajj(1,2)}{d(2,3)\akk}        &\frac{d(1,2)+d(2,3)}{d(1,3)}&\frac{d(1,4)}{d(1,3)}\\\hhline{|=*{2}{=:}=|}
\!\!\!\begin{tabular}{l}
Shortest path distance,\phantom{$1\ajj$}\\
   Resistance distance
      \end{tabular}\!\!\!            &                                    1&                           1&                    1.5\\\hhline{|=:*{3}{=|}}
Walk distance,             $\a=1\ajj$&                                 1.08&                           1&                                 1.52\\\hline
Long walk distance ($\a\to\infty$)   &\frac{1+\sqrt{5}\aj}{2\ak}\approx1.62&                           1&\frac{1+\sqrt{5}\aj}{2\ak}\approx1.62\\\hline
Logarithmic forest distance, $\a=  2$&                                 0.89&                       1\ajj&                   1.47\\\hhline{|=:*{3}{=|}}
Forest distance,           $\a=1\ajj$&                                 1.08&                        1.32&                                 1.26\\\hline
``Plain'' walk distance,     $\a=4.5$&                                 1.08&                        1.28&                                 0.95\\\hline
``Plain'' walk distance,     $\a=  1$&                                 0.96&                        1.46&                                 1.03\\\hline
\end{tabular}
$$
\Up{1.0}
\caption{The properties of several metrics on $P_4.$\label{ta_di}}
\end{table}

The walk distances and the logarithmic forest distances are graph-geodetic, so they satisfy $\frac{d(1,2)+d(2,3)}{d(1,3)}=1,$ since all paths between $1$ and $3$ visit~$2.$ Our examples suggest that these metrics are useful to model situations where, all other things being equal, the peripherality of vertices increases or decreases the distance between them. In such cases, the walk distances or the logarithmic forest distances can be used, respectively; in this example, for the former, ${d(1,2)>d(2,3),}$ while for the latter, ${d(1,2)<d(2,3).}$

The forest metrics \cite{CheSha97,CheSha00} are obtained by the application of \eqref{e_Da} to the matrices ${(I+\a L)^{-1},}$ where $L\!=\!\diag(A\bm1)\!-\!A$ is the Laplacian matrix of~$G.$ As well as walk metrics, they increase the distance between peripheral neighbors, however, the forest metrics are not graph-geodetic.

The ``plain'' walk metrics \cite{CheSha98} are obtained by the application of \eqref{e_Da} to the matrices $R_t={(I-t A)^{-1}}$, where $t=(\rho+\a^{-1})^{-1}$ (see\;\eqref{e_Rt} and\;\eqref{e_param}). Depending on $\a,$ they can either increase or decrease the distance between peripheral neighbors. Let us note that for $P_4$, they set $d(1,3)\approx d(1,4)$ or even $d(1,3)>d(1,4)$ (see the last column of Table\;\ref{ta_di}), which is quite exotic and does not meet the geodetic (graph traversal) approach taken in this paper.
%unnatural. So these metrics cannot be considered as a good model for applications.

Numerical examples and partial results suggest that the walk metrics and the logarithmic forest metrics more sensitively take into account the global structure of the graph than the electric metric\footnote{\emph{Electric metric\/} = resistance distance; this term was proposed by Harary (cf.\,\cite{KleinZhu98}).} does. In particular, the distances they provide depend not only on the paths between two vertices, but also on their centrality. As a result, these metrics do not coincide with the shortest path metric when $G$ is a tree.

\section*{Acknowledgements}

This work was partially supported by the RFBR Grant 09-07-00371 and the RAS Presidium Program ``Development of Network and Logical Control''.
The author is grateful to Ravindra Bapat, Michel Deza, and Ernesto Estrada for helpful discussions and to the anonymous referees for their comments.

\end{document}